\documentclass[a4paper, 11pt, reqno]{amsart}
\usepackage[protrusion=false,final]{microtype}
\usepackage{graphicx}
\usepackage{booktabs} 
\usepackage{pgfplots}
\usepackage{amssymb}
\usepackage{amsmath}
\usepackage{mathtools}
\usepackage{hyperref}
\usepackage{cleveref}
\numberwithin{equation}{section}
\usepackage{enumerate}
\usepackage{enumitem}
\usepackage{url}

\theoremstyle{plain}
\begingroup
\newtheorem{theorem}{Theorem}[section]
\newtheorem{lemma}[theorem]{Lemma}

\newtheorem{proposition}[theorem]{Proposition}
\newtheorem{corollary}[theorem]{Corollary}
\newtheorem*{thm*}{Theorem}
\endgroup
\theoremstyle{definition}
\begingroup
\newtheorem{definition}[theorem]{Definition}

\endgroup
\theoremstyle{remark}
\begingroup
\newtheorem{remark}[theorem]{Remark}
\endgroup

\renewcommand{\o}{\Omega}

\newcommand\e{\varepsilon}
\newcommand{\R}{\mathbb{R}}
\newcommand{\N}{\mathbb{N}}

\newcommand\iO{\int_\Omega}

\newcommand{\wk}[1]{\xrightharpoonup{#1}}

\newcommand{\liminfn}{\underset{n \to \infty}{\liminf}\;}
\newcommand{\limsupn}{\underset{n \to \infty}{\limsup}\;}
\newcommand{\limn}{\underset{n \to \infty}{\lim}\;}
\newcommand{\cur}{\mathrm{curl}\,}
\newcommand{\curd}{\mathrm{curl}_{\e}\,}
\newcommand{\curdn}{\mathrm{curl}_{\e_n}\,}
\newcommand{\dgrad}{\nabla_{\e}}
\newcommand{\dgradn}{\nabla_{\e_n}}
\newcommand{\Ly}{{L_{\e,\xi}}}
\newcommand{\Lyn}{{L_{\e_n,\xi}}}
\newcommand{\sm}{\setminus}
\newcommand{\cl}{\overline}
\newcommand{\M}{\R^{N\times N}}
\newcommand{\C}{\mathcal{C}}
\newcommand{\di}{{\rm div }\, }
\newcommand\gr{\nabla}
\newcommand{\curl}{\mbox{curl} \,}
\newcommand{\Div}{\mbox{div} \,}

\newcommand{\wkto}{\rightharpoonup}
\newcommand {\jump}[1] {[\![ #1 ]\!]}

\title{Vector Field Models for Nematic Disclinations}
\author[A. Acharya] {Amit Acharya} 
\address[A. Acharya]{Department of Civil and Environmental Engineering and Center for Nonlinear Analysis\\ Carnegie Mellon University\\
5000 Forbes Avenue\\Pittsburgh PA 15213, USA}
\author[I.Fonseca]{Irene Fonseca}
\address[I. Fonseca]{Department of Mathematics and Center for Nonlinear Analysis\\ Carnegie Mellon University\\
5000 Forbes Avenue\\Pittsburgh PA 15213, USA}
\author[L. Ganedi] {Likhit Ganedi} 
\address[L. Ganedi]{Institut für Mathematik\\ RWTH Aachen\\
Templergraben 55\\52062 Aachen, Germany}
\author[K.Stinson] {Kerrek Stinson} 
\address[K. Stinson]{Hausdorff Center for Mathematics\\ University of Bonn\\
Endenicher Allee 62 \\53115 Bonn, Germany}
\begin{document}
\maketitle
 \begin{abstract}
 In this paper, a model for defects that was introduced in \cite{ZANV} is studied. In the literature, the setting of most models for defects is the function space SBV (special bounded variation functions) (see, e.g., \cite{ContiGarroni, GoldmanSerfaty}). However, this model regularizes the director field to be in a Sobolev space by adding a second field to incorporate the defect. A relaxation result in the case of fixed parameters is proven along with some partial compactness results.
\end{abstract}

\section{Introduction}\label{sec:intro}

The purpose of this paper is to initiate the rigorous mathematical analysis of a model of the dynamics of disclination line defects in nematics proposed in \cite{ZANV}. Here, we focus on energetic aspects of the model. Combined with the ideas presented in \cite{AD} and the demonstrations provided in \cite{PAD, ZZAGW, ZANV}, which include static fields of straight $\pm \frac{1}{2}$  disclinations, their annihilation, the dissociation of closely bound pairs of straight disclinations, as well as static fields of disclination loops, the model can be considered as a thermodynamically consistent generalization of the Ericksen-Leslie (EL)  model to account for the dynamics of disclination lines, with total energy that remains bounded in finite bodies in the presence of these line defects.

The practical applications of liquid crystalline phases abound, from  liquid crystal displays for electronic devices and cell membranes in biology (use in the mechanically `soft' phase), to a vast variety of liquid crystal polymers including the mechanically `hard' Kevlar, for body armor, and in tires.  Equally, topological defects abound in liquid crystalline media, fundamentally due to microscopic structural symmetries related to the head-tail symmetry of the director. Due to this symmetry, a vector field assigned to a director field can undergo continuous changes in orientation around a non-unique surface terminating at a unique \textit{disclination line}, where the jump in orientation across the surface is quantized to be $\pi$ radians (see, e.g., De Gennes and Prost \cite[Sec.~4.2.1]{degennes-prost}). The line field representing the director and the corresponding vector field are both discontinuous at the disclination line, and the energy cost of such a line discontinuity can be sustained by the material. It is this fundamental insight, going back to the kinematics of singularities in linear elasticity theory due to Volterra and Weingarten (see, e.g., \cite{acharya2019weingarten} for a contemporary review), that forms the core idea of our model and, in fact, has been recently used to define an algorithm to detect line defects in molecular configurations of nematics produced by Molecular Dynamics simulations \cite{saha}, extending the notion of a director down to the level of a single nematic molecule. 

The understanding of the energetics and \emph{dynamics} of topological defects and their interaction form an important part of the theoretical study of liquid crystalline media, and we are particularly interested in the universality of this behavior across the material behaviors of liquid crystals and crystalline solids. A primary justification of the type of model we consider is that it is inter-operable with a study of dislocation defects in crystalline solids with simply a change in interpretation of the field variables involved (this cannot be said of the Landau-DeGennes model \cite{Sonnet-Virga,schopohl1987defect,svenvsek2002hydrodynamics} which, nevertheless, employs a more adapted order parameter, the $Q$ tensor, to describe the head-tail symmetry of the nematic director, albeit at the cost of including a biaxial phase within the description as well). A distinguishing feature of our model is that, even in a `smooth' setting, defect cores can be identified as a locally calculable field variable, arguably a desirable feature as discussed in \cite{biscari2005field}.

Our model introduces an extra second-order tensor field, $B$, beyond the EL director field, $k$ (strictly speaking $k$ is a vector field representation of the director line field). This new field is to be physically thought of as a locally integrable realization, at the mesoscale, of the `singular' part of the director gradient field ($Dk$) in the presence of line defects, singular when viewed at the macroscale. Thus, at the mesoscale, both the director gradient field and the new field are integrable
\footnote{The length scales $\xi$ (core width) and $\e \xi$ (the layer width) that appear subsequently in Sec.~\ref{sec:main_res} and are at the heart of such a physical regularization, can be precisely defined in configurations of nematic molecules arising in Molecular Dynamics simulations, as shown in ongoing work \cite{saha}.}
- with this clear, we nevertheless refer to $B$ as the `singular part of the director distortion.' Notably, the field $B$ is not a gradient, and this allows it to encode information on the topological charge of line defects through its curl. The evolution of the director field $k$ continues to be obtained from the balance of angular momentum, as shown by Leslie \cite{L1992}, and the evolution of $B$ follows from a conceptually simple conservation law for the topological charge of the line defects, which is tautological before the introduction of constitutive assumption for the disclination velocity, the latter deduced from consistency with the second law of thermodynamics. The introduction of dynamics based on such a conservation law, rooted in the kinematics of defect lines, is a conceptual departure from what is done for dynamics with the Landau-DeGennes $Q$ tensor model (see, \cite{Sonnet-Virga, macmillan_1}), or Ericksen's variable degree of orientation model \cite{Ericksen_1991}. In doing so, the model also makes connections to the dynamics of dislocation line defects in elastic solids \cite{ZAWB,AZA}, as well as their statics \cite{Arora_AA_2020}. 

At the length scales where individual line defects are resolved, partial differential equations-based dynamical models arising from continuum mechanical considerations involve Newtonian and thermodynamic driving forces that include nonlinear combinations of entities representing director distortions and the disclination density fields. This requires a minimum amount of regularity in these fields, and hence it is essential to have a formulation that utilizes at least locally integrable functions, and our model is designed to be consistent with this requirement (of course, this does not preclude the question of studying limiting situations of such models when such functions tend to singular limits, modeling fields that have discontinuities, and singularities in the macroscopic limit).

\subsection{Main Results}\label{sec:main_res}
We investigate the behavior of local minimizers for the previously discussed model for liquid crystals with disclination defects. Let $\Omega\subset \R^N$, where $N=2$ or $3$, be the domain occupied by a nematic liquid crystal. We consider the following energy for the director field $k \in W^{1,2}(\o;\R^2)$ and the singular part of the director distortion $B \in H_\text{curl}(\o;\R^{2\times 2})$,
 \begin{equation}\label{eq:energy}
 	E_{\e,\xi}[k,B] := \int_\o{\left[\frac{(|k|-1)^2}{\e\xi^2} + |\nabla k - B|^2 + \e \xi^2|\cur B|^2 + \frac{1}{\e \xi^2}W(\e \xi|B|)\right]dx}, 
 \end{equation}
where $W: [0,\infty) \to [0,\infty)$ is a nonconvex double-well continuous potential with wells at 0 and 2.
 
\begin{remark} \label{rmk:heuristics}
The heuristics behind the energy \eqref{eq:energy} above for the prediction $\pm \frac{1}{2}$ line defects are as follows: the nonconvex potential $W$ assigns vanishing energy cost when $|B| = 0$ or $ \frac{2}{\e \xi}$. This, along with the elastic energy term $| D k - B|^2$, assigns approximately vanishing elastic cost for pointwise values of the director gradient of the type $Dk \approx \frac{n - (-n)}{\e \xi} \otimes l$, where $0 < \e \leq 1$ and $n, l$ are unit vectors. Here, $n$ corresponds to the director field $k$, and $l$ represents the direction along which the jump of $n$ occurs. 

In Fig. \ref{fig:gcr}, this scenario is presented for a rectangular transition layer for the director field. 

 \begin{figure}
	\centering
	\includegraphics[width=0.5\linewidth]{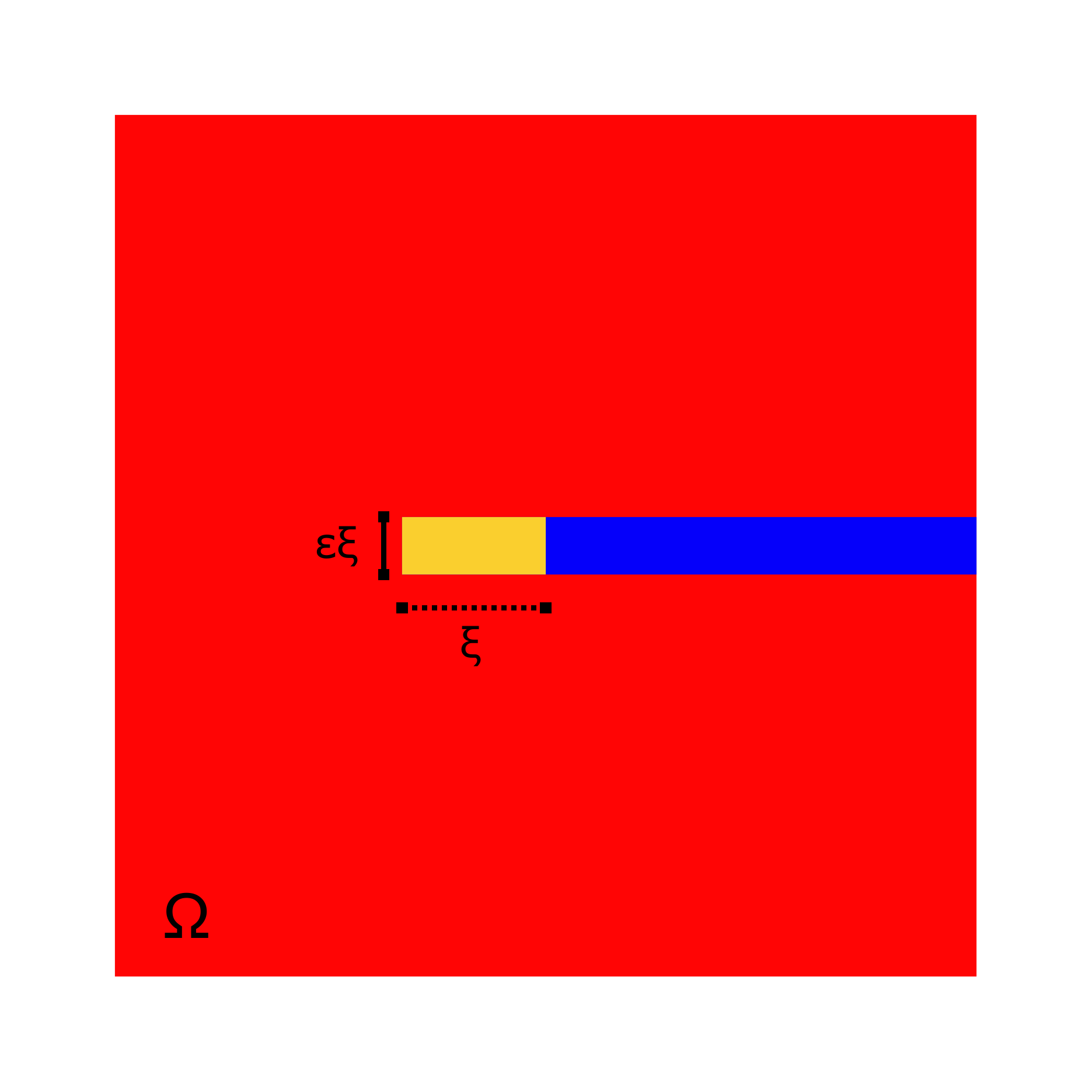}
	\caption{A representation of the domain and the layer where the discontinuity is supported}
	\label{fig:gcr}
\end{figure} 

If the transition layer in Fig. \ref{fig:gcr} did not terminate inside the domain, the energy cost would be minimal for a (diffuse) jump in the orientation of $k$ by $\pi$ radians across the layer. However, for a layer terminating inside the domain, $\curl B$ is non-zero near the termination (or core), and if the width of the layer in the vertical direction was $0$ (i.e., $\e = 0$), $\curl B$ would be singular (the classical defect solution results from the choice $\xi = 0$, when layer and core widths vanish). When the $\curl$ does not vanish, the density $Dk$ cannot annihilate $B$ (regardless of $\xi = 0$, or not). To see this, the Euler Lagrange equation of a functional with just the energy density $|Dk - B|^2$ for admissible variations in $k$ with $B$ a specified field is, with $Dk - B  =: e$, $\Div e = 0$, and $\curl e = -\curl B$. 

As an example, for $\curl B$ a (mollified) Dirac supported at the layer termination, this produces the approximate elastic energy density field, given here by $|e|^2$, of some canonical line defects in 2 dimensions (screw dislocation in solids, the wedge disclination in nematics with unit vector constraint imposed, either exactly or approximately) \cite{ZZAGW,ZANV, nabarro, hirth_lothe, frank1958liquid, degennes-prost}. Since $e = Dk$ outside the layer, we recover the relevant director field (using the penalized unit vector constraint represented by the first term in \eqref{eq:energy} and a specified value of $k$ at one point of the domain). Within the layer, but outside core, the director field $k$ flips orientation by $\pi$ radians, with a somewhat more involved distribution in the core.

\end{remark}
 
 As we are interested in minimizers of the energy (\ref{eq:energy}), we first consider the relaxation at a fixed scale (see Theorem \ref{thm:relax} for complete details). The energetic relaxation provides a functional to which the direct method of the calculus of variations is amenable, and is the first step to understanding the structure of minimizers.
 
 \begin{thm*}
Let $\Omega\subset \R^N$, $N=2$ or $3$, be an open, bounded set with Lipschitz boundary. For $\e,\xi>0$ fixed, the lower semicontinuous relaxation of the energy (\ref{eq:energy}) is given by
	\begin{align}
		\bar E_{\e,\xi}[k,B] &:= 
			\iO\left[\frac{(|k|-1)^2}{\e\xi^2}+|\nabla k - B|^2\right]\;dx \nonumber\\
			&\quad\quad + \iO\left[\e\xi^2|\cur B|^2 + \frac{1}{\e\xi^2}Q(W(|\cdot|))(\e\xi B)\right] \;dx.
	\end{align}
 for all $k \in W^{1,2}(\o;\R^N)$ and $B \in H_\cur(\o;\M)$.
	Here, $Qf$ denotes the quasiconvex envelope of $f$.
\end{thm*}

\begin{remark}	\label{rmk:quasiconvex}
	We show that $Q(W(|\cdot|))\equiv 0 $ in the ball $B(0,2)\subset \R^{N\times N}$. Clearly, $Q(W(|\cdot|))(p) \geq W^{**}(|p|)$ always. This shows that $Q(W(|\cdot|))(p) \geq 0$ everywhere. Further, we claim that each $p \in B(0,2)\subset \R^{N\times N}$ is given as the convex combination of two-elements in $2S^{N\times N-1}$ differing by a rank-one matrix. Indeed, note we can find $\lambda_+,\lambda_- \geq 0$ such that $$p_\pm:=p \pm \lambda_\pm e_1\otimes e_1 \in 2S^{N\times N-1}\quad\text{and}\quad \mathrm{rank}(p_+-p_-)=1.$$ 
		By taking $t= \frac{\lambda_-}{\lambda_++\lambda_-} \in [0,1]$, we can compute that 
	$$t(p + \lambda_+ e_1\otimes e_1) + (1-t)(p - \lambda_- e_1\otimes e_1) = p.$$
	As $W(|\cdot|)\llcorner_{2S^{N\times N-1}} \equiv 0$ by (\ref{hyp:f2}), we apply rank-one convexity of the quasiconvex envelope \cite{Dacoragna} and the fact that the envelope is always below the original function to find that 
 \begin{align*} 
 Q(W(|\cdot|)) (p) \leq & tQ(W(|\cdot|)) (p_+) + (1-t)Q(W(|\cdot|)) (p_-) \\
 \leq & tW(|p_+|) +(1-t)W(|p_-|) = 0,
 \end{align*} thereby concluding that $Q(W(|\cdot|))$ vanishes in the ball $B(0,2)$. 
 
 We further conjecture that $Q(W(|\cdot|))(p) = W^{**}(p)$ due to the radial symmetry, however characterization of the quasiconvex envelope poses challenges even in simple cases (see, e.g., \cite{QCenvelope} for one of the few nontrivial examples of a calculation of the envelope).
\end{remark}

To motivate the constraints we will place on the field $B$, we introduce a simple example. We now restrict our attention to dimension $N=2$ and consider the limit $\e \to 0$ with $\xi > 0$ fixed. Considering any $k \in C^{1}(\Omega;S^{1})$, we set $B := \nabla k$ to find that $\bar E_{\e,\xi}[k,\nabla k] = \int_\Omega\frac{1}{\e\xi^2}Q(W(|\cdot|))(\e\xi \nabla k) \,dx.$ Given Remark \ref{rmk:quasiconvex} above, the function $Q(W(|\cdot|))(\e\xi B)$ doesn't see the energy from $B$ if $|B|\leq \frac{2}{\e \xi}$. As a result, if $\epsilon \leq \frac{2}{\xi\|\nabla k\|_\infty +1}$, then $\bar E_{\e,\xi}[k,\nabla k]=0$. Such a result (though defect free) shows that further constraints on the field $B$ are required to gain physically meaningful insight in the limit as $\e\to 0$.

We consider the particular case of $B$ supported in a layer as in Figure \ref{fig:gcr}, a physically relevant geometric configuration (see, e.g., \cite{ZANV}). 

To be precise, let $\Omega:= (-1,1)^2$ be the domain of a liquid crystal in the plane. We assume the defect is at the origin, and the surface of discontinuity is within a layer $L_{\varepsilon,\xi}$, with 
\begin{equation}\label{eqn:BlayerConstraint}
\Ly:= (-\xi,1)\times\left(-\frac{\e\xi}{2},\frac{\e\xi}{2}\right) \quad \text{ and } \quad B \equiv 0 \text{ in }\Omega \setminus L_{\e,\xi} 
\end{equation} with parameters $\e,\xi > 0$. In physical terms, $\xi$ is the core length of the crystalline defect and $\e$ is a parameter determining the thickness of the defect layer $L_{\e,\xi}$ (see Figure \ref{fig:gcr}).  In this paper, we are primarily concerned with $\pm\frac{1}{2}$ disclinations, which must satisfy the constraint
  \begin{equation}\label{eq:halfdisclin}
 \left|\int_\o{\cur B\;dx} \right| = 2.
 \end{equation}
 This is a model constraint requiring a disclination to exist in the domain. By Stokes' theorem, (\ref{eq:halfdisclin}) is consistent with a layer field in the form of $ B = \frac{n - (-n)}{\e \xi} \otimes l$ in a layer of width $\e \xi$ with normal in the direction $l$ and $n$ a unit vector corresponding to the director field $k$ (see Fig. \ref{fig:gcr}), and as described in Remark \ref{rmk:heuristics}.

After a change of variables analogous to typical dimension reduction problems \cite{DretRaoult}, we prove a compactness theorem for the rescaled fields, which is precisely stated in Theorem \ref{thm:compactness}. Next, we state a theorem that follows from Theorem \ref{thm:compactness} and which emphasizes the coupling of the physical quantities in the asymptotic limit.

\begin{theorem}\label{thm:compactCorollary}
    Let $k_\e$ and $B_\e$ have uniformly bounded energy as $\e \to 0$, that is, $\sup_{\e>0}E_{\e,\xi}[k_{\e},B_{\e}]<\infty$. Further, suppose $B_{\e}$ satisfies the geometric constraint (\ref{eqn:BlayerConstraint}) and corresponds to $\pm 1/2$ disclination by satisfying (\ref{eq:halfdisclin}). Then up to a subsequence (not relabeled), $k_\e \to k$ strongly in $L^2(\o;\R^2)$ where $k \in SBV(\o;S^1)\cap W^{1,1}(\Omega\setminus L_{0,\xi};S^1)$. Defining the jump of $k$ on $(-\xi,1)$ by $\jump{k}: = k^+ -k^-$, the compatibility condition
	$$\jump{k}(s) = \int_{-\xi}^{s}\int_{-\frac{\xi}{2}}^{\frac{\xi}{2}} \alpha \;dx \quad \text{ and }\quad |\jump{k}(1)| = 2 $$
	is satisfied for all $s\in (-\xi,1)$, where $\alpha\in L^2(L_{1,\xi})$ is the limit of the rescaled $\{\tilde B_\e := \e B(x_1,\e x_2)\}_{\e}.$ Furthermore, $\nabla k$, the part of $Dk$ absolutely continuous with respect to the Lebesgue measure has higher regularity, in the sense that
    $$\cur(\nabla k) = - \frac{d}{dx_1}\jump{k}\otimes e_2\mathcal{H}^{1}\llcorner_{(-\xi,1)\times\{0\}}.$$ 
\end{theorem}

We can also make a connection to the recent preprint \cite{GoldmanSerfaty}, where a SBV model for $\pm \frac{1}{2}$ disclinations is proposed and the constraint that $[k] = 2$ along the jump set is imposed. The energy we use can be viewed as an attempt to also relax the one used in \cite{GoldmanSerfaty} by being a Sobolev model allowing for a more general class of jumps in the SBV limit.

There are many open questions stemming from this work which we highlight in Section \ref{sec:conclusion}. Foremost, is the integral representation of a precise limiting energy for the case $\e \to 0$ with $\xi > 0$ fixed. It is also possible to consider the case of $\xi \to 0$ at various rates compared to $\e \to 0$. However, the limit $\xi \to 0$ will be complicated by the need to rescale the energy by $\log \xi$, which leads to a delicate Ginzbug-Landau type problem (see \cite{JerrardSoner}, \cite{AlicandroPonsig}).  

\section{Mathematical Preliminaries}

Let $W^{1,2}(\o;\R^2)$ denote the usual Sobolev space, and we designate by $H_\text{curl}(\o;\R^{2\times 2})$ the space of $L^2$ matrix valued tensors, whose row-wise distributional curl is also in $L^2$. Under this setting, we consider the energy (\ref{eq:energy}) with a nonconvex continuous potential
$W: [0,\infty) \to [0,\infty)$ satisfying the following coercivity and growth properties,
\begin{equation}\label{hyp:f1}
	\frac{1}{C}|x|^2-C\leq W(x) \leq C(1+|x|^2) \quad \text{for }\; x \in [0,\infty)
\end{equation}
for some $C>0$, and
\begin{equation}\label{hyp:f2}
\{x \in [0,+\infty): W(x) = 0\} = \{0,2\}.
\end{equation}

The mathematical framework we use to study the convergence of the functional \eqref{eq:energy} is encapsulated by the notion of $\Gamma$-convergence, which we recall next. 

\begin{definition}\label{def:Gamma_conv}
Given a metric space $(X,d)$, let $F_n: X \to [0,\infty]$ be a sequence of functionals for $n\in \N$. We say that $F_n$ $\Gamma$-converge to $F_0: X \to [0,\infty]$ with respect to the metric $d$ if the following two conditions hold:
\begin{enumerate}
    \item (Liminf Inequality) For every $u \in X $ and for every sequence $\{u_n\} \subset X $ such that $u_n \to u$ with respect to the metric $d$, we have $$F_0(u) \leq \liminfn{}{F_n(u_n)}.$$
    \item (Recovery Sequence) For every $u \in X$, there exists $\{u_n\} \subset X$ such that $u_n \to u$ with respect to the metric $d$, and the sequence recovers the energy, i.e., $$\limsupn{F_n(u_n)}= F_0(u).$$
\end{enumerate}
\end{definition}

For the relaxation of the energy (\ref{eq:energy}), which is by definition the $\Gamma$-limit of the constant sequence of functionals $F_n :=E_{\e,\xi}$, we will rely on the now classical notion of quasiconvexity introduced by Morrey \cite{Morrey}. Analogous to the characterization of convex functions via Jensen's inequality, quasiconvex functions satisfy a Jensen's type inequality for gradient fields. Specifically, a Borel measurable function $g: \R^{d\times N}\to [0,\infty]$, $d,N\geq 1$, is quasiconvex if
$$g(\xi) \leq \int_{(0,1)^N} g(\xi + \nabla \phi) \, dx \quad \quad \text{ for all }\xi \in \R^{d\times N} \text{ and }\phi \in C^1_c((0,1)^N;\R^d).$$
The quasiconvex envelope of a function $f$ is given by the greatest quasiconvex function beneath $f$, i.e., it is defined pointwise by
\begin{equation}\label{def:quasiconEnv}
Qf(x_0) : = \sup\{g(x_0):g \text{ is quasiconvex and }g\leq f\}.
\end{equation}
We refer the reader to \cite{Dacoragna} for further details on such functions.

As mentioned in the discussion preceding Theorem \ref{thm:compactCorollary}, we wish to model discontinuities across 2-d surfaces - when viewed at the macroscale -  in a vector field representation of the director field containing a line defect, while incorporating the fact that at the microscopic scale such a jump across the surface must necessarily be spread over a region roughly of the order of the spacing between adjacent mesogens (cf. \cite{saha}). Specifically, we impose the condition that $B$ vanishes outside of the layer, which is equivalent to the conditions
   \begin{align}
  	&B = 0 \quad\quad \mathrm{ in }\quad \o \setminus L_{\e,\xi}, \label{eq:B_zero}\\
  	&Bt = 0\quad\quad \mathrm{ on }\quad \partial \Ly \sm \partial \o, \label{eq:B_tan}
  \end{align}
where $t$ is the tangent vector to the boundary point. The condition \eqref{eq:B_tan} comes for free as $B \in H_{\cur}(\o;\R^{2\times 2})$ and there is a well-defined tangential trace matching the condition \eqref{eq:B_zero} \cite{FabrieBoyer}.  

Finally, we recall that a function $u$ belongs to $BV(\Omega)$ if its distributional gradient is given by a finite Radon measure. Informally, in the case that $u$ has only surface discontinuities, $u$ belongs to $SBV(\Omega)$, and for the sake Theorem \ref{thm:compactCorollary}, it suffices to know that if $u \in W^{1,1}(\Omega \setminus K)\cap L^\infty(\Omega)$, where $K$ is a closed set with finite surface measure, i.e., $\mathcal{H}^{N-1}(K)<\infty,$ then $u \in SBV(\Omega)$. For further details, we refer to \cite[Proposition 4.4]{AmbrosioFuscoPallara}.

\section{Relaxation for fixed $\e,\xi$}
We study properties of the energy (\ref{eq:energy}) with $\e,\xi>0$ fixed. A priori, it is not clear that minimizers to the problem exist nor is it clear what the value of the infimum is. In order to apply the direct method of the calculus of variations, we must consider the lower semicontinuous envelope of the functional. Specifically, we obtain an integral representation for the relaxation of energy (\ref{eq:energy}) in dimensions $N=2$ or $3$. This dimension constraint enables us to use the Helmholtz decomposition and the corresponding Sobolev spaces, as detailed in \cite{FabrieBoyer}. Here, $Qf$ denotes the quasiconvex envelope of $f$ as in (\ref{def:quasiconEnv}).
\begin{theorem}\label{thm:relax}
	Let $\Omega\subset \R^N$, with $N=2$ or $3$, be an open, bounded set with Lipschitz continuous boundary, and let $E$ be defined in \eqref{eq:energy}. For all $k \in W^{1,2}(\o;\R^N)$ and $B \in H_\cur(\o;\M)$, the relaxation of $E_{\e,\xi}$ is given by
	\begin{equation}\nonumber
		\bar E_{\e,\xi} [k,B] : = \inf \left\{\liminfn E_{\e,\xi}[k_n, B_n]: (k_n,B_n) \to (k,B) \right\},
	\end{equation}
	where  the convergence is such that $k_n \to k$ and $B_n \wk{} B$ in $L^2$. 
	This relaxed energy has the integral representation
	\begin{align}\label{def:Ebar}
		\bar E_{\e,\xi}[k,B] &:= 
			\iO\left[\frac{(|k|-1)^2}{\e\xi^2}+|\nabla k - B|^2\right]\;dx \nonumber\\
			&\quad\quad + \iO\left[\e\xi^2|\cur B|^2 + \frac{1}{\e\xi^2}Q(W(|\cdot|))(\e\xi B)\right] \;dx.
	\end{align}
	
\end{theorem}

We note that this result can equivalently be phrased as $E_{\e,\xi}$ $\Gamma$-converge to $\bar E_{\e,\xi}.$
To prove Theorem \ref{thm:relax}, we will introduce an intermediate functional, related to (\ref{def:Ebar}) through the Helmholtz decomposition of $B$.
We denote the space of the divergence-free fields with integrable $\cur$ as
\begin{equation}\label{def:Cspace}
	\C : = \{u \in H_{\cur}(\Omega;\M) \; : \ \di u = 0, \ u \cdot \nu = 0 \ \text{on } \partial \Omega \}.
\end{equation}
Define the functional $I: (L^2(\Omega;\R^N))^2 \times L^2(\Omega;\M) \to [0,+\infty]$ by
\begin{equation}\label{def:I}
\begin{aligned}
    	I[\tilde{k},z,p] := &\int_\Omega\left[(|\tilde{k}+z|-1)^2+|\nabla \tilde{k}|^2\right]\;dx \\
		&\quad +\iO\left[|p|^2+|\cur p|^2 + W(|\nabla z + p|)\right]\;  dx,
\end{aligned}
\end{equation}
if $(\tilde{k},z,p) \in \mathcal{X}$, and $+\infty$ otherwise.
Here,
\begin{equation}\label{def:Xspace}
	\mathcal{X}:= W^{1,2}(\Omega;\R^N)\times \left(W^{1,2}(\Omega;\R^N) \cap \left\{\int_\Omega z \, dx = 0\right\}\right) \times \C.
\end{equation}

First, we investigate compactness of the functional $I$ in the following lemma.
\begin{lemma}[Compactness of $I$]\label{lem:compact}
	Let $\Omega \subset \R^N$, with $N=2$ or $3$, be an open, bounded set with Lipschitz continuous boundary. Consider a sequence $\{(\tilde{k}_n,z_n,p_n)\}$ such that $\sup_{n \in \N }{I[\tilde{k}_n,z_n,p_n]}\leq C$.
	Then there is $(\tilde{k},z,p) \in \mathcal{X} $ such that, up to a subsequence (not relabeled),
	\begin{align*}
		\tilde{k}_n &\wk{} \tilde{k} \quad \text{in } W^{1,2}(\Omega;\R^N),\\
		z_n &\wk{} z \quad \text{in } W^{1,2}(\Omega;\R^N),\\
		p_n &\wk{} p \quad \text{in } W^{1,2}(\Omega;\M),
	\end{align*}
	and
	$$(\tilde{k}_n,z_n,p_n) \to (\tilde{k},z,p) \quad \text{strongly in } (L^2(\Omega;\R^N))^2 \times L^2(\Omega;\M).$$
\end{lemma}
\begin{proof}
	The key is the inequality (see \cite{FabrieBoyer})
	\begin{equation}\label{bdd:div+curl}
		\|p_n\|_{W^{1,2}(\Omega; \R^N)} \leq C(\Omega)\left(\|p_n\|_{L^2(\Omega;\R^N)} + \|{\rm div }\, p_n \|_{L^2(\Omega)} + \|\cur p_n \|_{L^2(\Omega;\R^N)} \right).
	\end{equation}
	From the definition of $\C$ in (\ref{def:Cspace}), we have that $\di p_n = 0$. As $I$ in (\ref{def:I}) controls $\|\cur p_n\|^2_{L^2(\Omega)}$, we apply the above inequality to conclude that $\|p_n\|_{W^{1,2}(\Omega,\M)}\leq C<\infty$. Convergence as in the statement of the lemma follows from weak compactness. By (\ref{hyp:f1}) and (\ref{def:I}), we have $\sup_n\|\nabla z_n\|_{L^2(\Omega)} \leq C<\infty,$ and the desired convergence follows from Poincar\'e's inequality because $\int_\Omega z_n \, dx = 0.$ Finally, the uniform bound of the energy (\ref{def:I}) implies a uniform bound on $\|\nabla \tilde k_n\|_{L^2(\Omega)}$. Combining this with control of $\|z_n\|_{L^2(\Omega;\R^N)}$ and $\|\tilde k_n\|_{L^2(\Omega;\R^N)}$, shows that $\sup_n\|\tilde k_n\|_{W^{1,2}(\Omega;\R^N)}\leq C<\infty$.
	
	To conclude strong convergence in $L^2,$ we apply the Rellich-Kondrachov compactness theorem.
\end{proof}

We will prove the relaxation of the functionals $I$ using techniques for \\
$\Gamma-$convergence (see Definition \ref{def:Gamma_conv}), i.e., that $I$ $\Gamma$-converges to $\bar{I}$, with
\begin{equation}\label{def:I_limit}
\begin{aligned}
    	\bar I[\tilde{k},z,p] :=& \int_\Omega \left[(|\tilde{k}+z|-1)^2+|\nabla \tilde{k}|^2\right]\;dx \\
		& \quad +\iO\left[|p|^2+|\cur p|^2 + Q(W(|\cdot|))( \nabla z + p)\right]\;  dx,
\end{aligned}
\end{equation}
if $(\tilde{k},z,p) \in \mathcal{X}$, and $+\infty$ otherwise.

We note that this functional is the same as the original functional with $W$ replaced by $Q(W(|\cdot|))$. In order to prove that this is indeed the correct limiting energy, we first show that the $\liminf$ inequality of $\Gamma-$convergence is satisfied (see Definition \ref{def:Gamma_conv}(1)).

 \begin{lemma}[Liminf of $I$]\label{lem:liminfI}
 	Let $\Omega\subset \R^N$, with $N=2$ or $3$, be an open, bounded set with Lipschitz continuous boundary, and assume that $W$ satisfies (\ref{hyp:f1}) and (\ref{hyp:f2}). For all sequences such that $$(\tilde{k}_n,z_n,p_n)\to (\tilde{k},z,p) \quad \text{strongly in } (L^2(\Omega;\R^N))^2 \times L^2(\Omega;\M),$$ we have
 	$$
 	\bar{I}[\tilde{k},z,p]\leq \liminfn I[\tilde{k}_n,z_n,p_n].
 	$$
 	
 \end{lemma}
 
 \begin{proof}
 	We define the function $h(s,\eta) : = W(|s + \eta|)$ and $Q h(s,\cdot)$ to be the greatest quasiconvex function below $h(s,\cdot)$ as in (\ref{def:quasiconEnv}). 
  
    We claim that
 	\begin{multline}\label{bdd:hLSC}
 		\int_\Omega Q(W(|\cdot|))(p+\nabla z ) \, dx = \int_\Omega Q h(p,\nabla z) \, dx \\
 		\leq \liminfn \int_\Omega h(p_n,\nabla z_n) \, dx = \liminfn \int_\Omega W(|p_n +\nabla z_n|)\, dx.
 	\end{multline}
 	The first equality is easy to see because we note that there is a translational symmetry to $h$ which gives the equality 
 	\begin{equation}\label{eqn:quasiTranslation}
 	Q h(s,\nabla z) = Q h(0,\nabla (sx+z(x))) = Q(W(|\cdot|))(s+\nabla z )
 	\end{equation}
 	for all $s \in \R$ and almost every $x \in \R^N$; the first inequality in (\ref{bdd:hLSC}) follows choosing $s = p(x)$ pointwise. To obtain the second inequality, note that $Qh(s,\eta)$ is $2$-Lipschitz in the second variable, so by (\ref{eqn:quasiTranslation}), it is continuous in the $s$-variable and hence $Qh$ is Carath\'eodory. Applying the lower semicontinuity result of Acerbi and Fusco for quasiconvex functions \cite{AcerbiFusco}, we find that 
 	\begin{equation}\nonumber
 	\int_\Omega Q h(p,\nabla z) \, dx \leq \liminfn \int_\Omega Qh(p_n,\nabla z_n) \, dx	\leq \liminfn \int_\Omega h(p_n,\nabla z_n) \, dx,
 	\end{equation}
 	concluding (\ref{bdd:hLSC}).

 The lower semi-continuity of the $L^2$ norm shows that
 	\begin{equation}\label{bdd:LSCL2}
 		\liminfn \int_\Omega \left(|\nabla \tilde k_n|^2+|\cur p_n|^2 \right) \, dx \geq \int_\Omega \left(|\nabla \tilde k|^2 + |\cur p|^2\right) \, dx,
 	\end{equation}
 	with strong convergence of $z_n$, $k_n$, and $p_n$ in $L^2$ giving that
 	\begin{equation}\label{eqn:strConv}
 		\lim_{n \to \infty} \int_\Omega \left[(|\tilde{k}_n + z_n|-1)^2  + |p_n|^2\right] dx = \int_\Omega \left[(|\tilde{k}+z|-1)^2+ |p|^2 \right] dx.
 	\end{equation}
 	Combining (\ref{bdd:hLSC}), (\ref{bdd:LSCL2}), and (\ref{eqn:strConv}), the lemma is proven.
 \end{proof}

In order to complete the integral representation for the relaxation of $I$, we show the existence of a recovery sequence (see Def \ref{def:Gamma_conv}(2)).

\begin{lemma}[Recovery Sequence for $I$]\label{lem:limsupI}
	If $(\tilde{k},z,p) \in \mathcal{X}$, then there exists a sequence $\{(\tilde{k}_n,z_n,p_n)\}$ such that
	\begin{align}
		(\tilde{k}_n,z_n,p_n)\to (\tilde{k},z,p) \quad &\text{strongly in } (L^2(\Omega;\R^N))^2 \times L^2(\Omega;\M)\label{cv}\\
		\limsupn I[\tilde{k}_n,&z_n,p_n]\leq \bar{I}[\tilde{k},z,p]. \label{ls}
	\end{align}
	
\end{lemma}

\begin{proof}
    As $\tilde{k}$ and $p$ are admissible in the original energy, we take $p_n := p$ and $\tilde{k}_n := \tilde{k}$. 
	Now define the function $g:\Omega \times \M \to [0,\infty)$ by $g(x,\eta) = W(|p(x)+\eta|)$. Note that as $W$ is continuous and $p$ is measurable, $g$ is a Carath\'eodory function and has polynomial growth in $\eta$.
	By standard relaxation results \cite{Dacoragna}, we can find a sequence $ \{z_n\}\subset W^{1,2}(\Omega;\R^N)$ such that
	\begin{align}
		&z_n \to z \quad \text{strongly in } L^2(\Omega;\R^N), \label{str}\\
		&\limsupn \iO{g(x,\gr z_n(x))\; dx}\leq \iO{Qg(x,\gr z(x))\; dx}\label{qc}
	\end{align}
    where $Qg(x,\eta)$ is defined in (\ref{def:quasiconEnv}). Again, by a similar argument for (\ref{eqn:quasiTranslation}), we have that $$Qg(x,\gr z(x)) =  Q(W(|\cdot|))(p(x)+\nabla z(x)) $$ for $x\in \Omega$ almost everywhere.
	By the above relation, \eqref{str}, and \eqref{qc}, we obtain (\ref{ls}) as the other functions in $I$ (\ref{def:I}) are fixed.
\end{proof}

Combining the last two lemmas, we finish the proof of Theorem \ref{thm:relax}.

\begin{proof}[Proof of Theorem \ref{thm:relax}]
As $\e$ and $\xi$ are fixed, we will let $\e = \xi = 1$ without loss of generality.
Note that for any $B \in L^2(\Omega; \R^{N \times N})$ we can apply the Helmholtz decomposition (see \cite{FabrieBoyer}) row-wise to find $B = p + \nabla z$ for some $p \in \mathcal{C}$ that is divergence free (see (\ref{def:Cspace})) and $z \in W^{1,2}(\Omega;\R^N).$ As the $\cur$ vanishes on gradients and $p\in \mathcal{C}$ and $\nabla (z-k)$ are orthogonal in $L^2(\Omega;\M)$ by an integration by parts as $\di  p = 0$ and $p\cdot \nu = 0$ (see (\ref{def:Cspace})), we have 
\begin{align*}
E_{1,1}[k,B] &= E_{1,1}[k,p+\nabla z]  \\
& = \int_\Omega \left[(|k|-1)^2 + |\nabla k - (p+\nabla z)|^2 + |\cur p |^2 + W(p + \nabla z)\right] dx\\
&= I[k-z,z,p],
\end{align*}
where we have used the fact that
\begin{align*} \int_\Omega |\nabla k - (p+\nabla z)|^2 \ , dx &= \int_\Omega \left[|\nabla k - \nabla z|^2 + 2 p \cdot\nabla( z-  k)   + |p|^2 \right] dx \\
&= \int_\Omega \left[|\nabla k - \nabla z|^2  + |p|^2\right] dx,
\end{align*}
and likewise
$$\bar E_{1,1}[k,B] = \bar I[k-z,z,p].$$
Using the above relations, one has that Lemma \ref{lem:liminfI} and Lemma \ref{lem:limsupI} translate to $\liminf$ and $\limsup$ relations for $E_{1,1}$, thereby proving $E_{1,1}$ $\Gamma$-converges to $\bar E_{1,1}$ and concluding the proof.
\end{proof}

\section{Constrained Minimizers}
In the following, we study low energy sequences for the energy $E_{\e,\xi}$ in dimension $N=2$ as $\e \to 0$, with $\xi>0$ fixed, when the fields $B$ satisfy geometric (\ref{eqn:BlayerConstraint}) and defect (\ref{eq:halfdisclin}) constraints. 

Here the disclination layer $\Ly$ becomes thin in the limit. As it is typical in dimension reduction problems, we perform the change of variables
\begin{align}\label{eq:Change_of_Var}
	\tilde{k}(x_1,x_2):=k\left(x_1,\e x_2\right),\\
	\tilde{B}(x_1,x_2):= \e B\left(x_1,\e x_2\right).
\end{align}
We remark that we have rescaled the field $B$ as well because the quadratic coercivity of $W$ only gives control over $\e B$. 
With $L_\xi:= (-\xi,1)\times(-\frac{\xi}{2},\frac{\xi}{2})$, we write the energy in terms of a bulk contribution and layer contribution as

 \begin{align}
	&E_{\e,\xi}[k,
	B] = \int_{\o\sm\Ly}{\left[\frac{(|k|-1)^2}{\e\xi^2} + |\nabla k|^2 \right]\;dx}  \nonumber\\
	& + \e \int_{L_\xi} {\left[\frac{(|\tilde{k}|-1)^2}{\e\xi^2} + \left|\dgrad \tilde{k} - \frac{\tilde{B}}{\e}\right|^2+\frac{\xi^2}{\e}|\curd\tilde{B}|^2 + \frac{1}{\e \xi^2}W(\xi|\tilde{B}|)\right]\; dx} \nonumber\\
	&= \int_{\o\sm\Ly}{\left[\frac{(|k|-1)^2}{\e\xi^2} + |\nabla k|^2 \right]\;dx} \nonumber\\
	&\quad + \int_{L_{\xi}} {\left[\frac{(|\tilde{k}|-1)^2}{\xi^2} + \e\left|\dgrad \tilde{k} - \frac{\tilde{B}}{\e}\right|^2+\xi^2|\curd \tilde{B}|^2 + \frac{1}{\xi^2}W(\xi|\tilde{B}|)\right]\; dx} \nonumber\\
	 & =: E^{bulk}_{\e,\xi}[k] + E^{layer}_{\e,\xi}[\tilde{k},\tilde{B}] \label{eq:bulk+layer}
\end{align}
where $\dgrad := [\partial_1,\frac{1}{\e}\partial_2]$ and the scaled curl operator is $\curd g := \partial_1 g_2 - \frac{1}{\e}\partial_2 g_1$.
Furthermore, the $\cur$ constraint \eqref{eq:halfdisclin} becomes
 \begin{equation}\label{eq:curl_contraint_cov}
	\left|\int_{L_\xi}{\curd \tilde{B}\;dx} \right| = 2.
\end{equation}
Finally, we define the limit layer 
\begin{equation}\label{eqn:limitLayer}
L^0_\xi : = (-\xi,1)\times \{0\}.
\end{equation}

In the following, we take an arbitrary subsequence $\e_n \to 0$ and investigate sequences with uniformly bounded energy with an eye towards ultimately understanding the most physically relevant effective energies
\begin{equation} \label{eq:limit}
	E_\xi [k] : = \inf{\{\liminfn {E_{\e_n,\xi}[k_n,\tilde{B}_n]}: k_n \chi_{(L_{\e_n,\xi})^c} \to k \quad \mathrm{strongly \;in}\;\; L^2(\o;\R^2)\}}
\end{equation} 
and
\begin{align*} \label{eq:limit}
	E_\xi [k,\alpha] : = \inf{\{\liminfn {E_{\e_n,\xi}[k_n,\tilde{B}_n]}: k_n \chi_{(L_{\e_n,\xi})^c} \overset{L^2}{\to} k, \quad \curdn \tilde{B}_n \wk{L^2} \alpha \}}.
\end{align*} 

\subsection{Compactness}

Under the hypotheses (\ref{eqn:BlayerConstraint}) and (\ref{eq:halfdisclin}), we consider any sequence with uniformly bounded energy, and using (\ref{eq:Change_of_Var}) and (\ref{eq:bulk+layer}), we write it as the sum of the non-negative energies 
\begin{equation}\label{eq:unifBddEnergy}
 E^{bulk}_{\e_n,\xi}[k_n] + E^{layer}_{\e_n,\xi}[\tilde{k}_n,\tilde{B}_n] = E_{\e_n,\xi}[k_n, B_n] \leq C<\infty.
 \end{equation}

We prove the following theorem. 
\begin{theorem}\label{thm:compactness}
    Let $k_n,\tilde{B_n}$ have uniformly bounded energy  as in (\ref{eq:unifBddEnergy}). Then the following hold:
    \begin{itemize}
        \item Outside the layer, $k_n \chi_{(L_{\e_n,\xi})^c} \to k$ strongly in $L^2(\o;\R^2)$, where $k \in SBV(\o;S^1)$. Furthermore, the jump set of $k$ is contained in $\cl{L^0_\xi}.$
        \item In the layer, we can generate a rescaled $k_n$ which is denoted by $\tilde{k}_n$. 
    We will have the convergences:
    \begin{align*}
    &\tilde{k}_n \wk{} \tilde{k} \quad \mathrm{weakly\;\; in\;}  L^2(L_\xi;\R^{2}), \\
    &\tilde{B}_n \wk{} B \quad  \mathrm{weakly\;\; in\;} L^2(L_\xi;\R^{2 \times 2}),\\
    &\curdn \tilde{B}_n \wk{} \alpha \quad \mathrm{weakly\;\; in\;}  L^2(L_\xi;\R^{2}),
    \end{align*}
    for some $\tilde{k} \in  L^2(L_\xi;\R^{2})$, $\alpha \in L^2(L_\xi;\R^{2})$, and $B := \begin{bmatrix}
	0 & \partial_2 \tilde{k}_1 \\
	0 & \partial_2 \tilde{k}_2
\end{bmatrix}.$ 
\item Further, for almost every $x_1 \in (-\xi,1)$, $\tilde k(x_1,\cdot) \in W^{1,2}(-\frac{\xi}{2},\frac{\xi}{2})$ and we may define $[\tilde{k}]: (-\xi,1) \to {\R^2}$ as 
\begin{equation}\label{def:tildekjump}
[\tilde{k}](x_1) := \int_{-\frac{\xi}{2}}^{\frac{\xi}{2}}{\partial_2\tilde{k}\;dx_2}.
\end{equation} The following compatibility condition between $\tilde{k}$ and $\alpha$ holds:  
	$$[\tilde{k}](s) = \int_{-\xi}^{s}\int_{-\frac{\xi}{2}}^{\frac{\xi}{2}} \alpha \;dx,\qquad\qquad |[\tilde{k}](1)| = 2. $$
    \end{itemize}
\end{theorem}

We note that even though $k$ and $\tilde{B}$ were independent in the beginning, these fields become coupled through $\tilde k$ in the limit. In Proposition \ref{prop:traceCoincidence}, we will show that this relation can be directly expressed without the intermediate field $\tilde{k}$.

\begin{proof}
\textbf{Step 1: Bulk Energy.}
To control the energy in the bulk, we observe that 
\begin{equation}\label{eq:Uniform_bds_Bulk}
	\sup_n \int_{\o\sm\Lyn}{\left[\frac{(|k_n|-1)^2}{\e_n\xi^2} + |\nabla k_n|^2 \right]\; dx} \leq C .
\end{equation}
In particular, for any $U$ smooth open set which is compactly contained in the set $\o\sm\cl{L^0_\xi}$ (recall (\ref{eqn:limitLayer})), we have that $$\sup_n\|k_n\|_ {W^{1,2}(U;\R^2)}\leq C.$$
Thus, up to a subsequence (not relabeled), we have that $k_n \to k$ strongly in $L^2(\o\sm\cl{L^0_\xi};\R^2)$.
Because of the unit norm regularization, we have that $|k|=1$ almost everywhere. Furthermore, since $k \in W^{1,2}(\o \sm\cl{L^0_\xi};\R^2)$, it is an integration by parts argument \cite[Proposition 4.4]{AmbrosioFuscoPallara} to show $k \in SBV(\o; \R^2)$ where the jump set of $k$ is contained in $\cl{L^0_\xi}$ up to a set of $\mathcal{H}^1$-measure zero.

\textbf{Step 2: Layer Energy.}
In this portion of the energy, we have
\begin{multline}\label{eq:Uniform_bds_Layer}
    	\int_{L_\xi} \left[\frac{(|\tilde{k}_n|-1)^2}{\xi^2}  + \e_n\left|\nabla_{\e_n} \tilde{k}_n - \frac{\tilde{B}_n}{\e_n}\right|^2\right]\;dx \\
    	+\int_{L_\xi}\left[\xi^2|\cur_{\e_n}\tilde{B}_n|^2 + \frac{1}{\xi^2}W(\xi |\tilde{B}_n|)\right]\; dx \leq C.	
\end{multline}

 Using the quadratic coercivity of $W$ in \eqref{hyp:f1}, we have 
\begin{align}\label{eq:bounds_layer}
	&\|\tilde{k}_n\|_{L^2} \leq C, \\
	&\|\e_n \dgradn \tilde{k}_n - \tilde{B}_n  \|_{L^2} \leq C \e_n^\frac{1}{2} \label{eq:k_b_relation},\\
	&\|\curdn\tilde{B}_n\|_{L^2} + \|\tilde{B}_n\|_{L^2} \leq C. \label{eq:B_bounds}
\end{align}
This implies that up to a subsequence, not relabeled, we obtain 
\begin{align}
&\tilde{B}_n \wk{} B \quad  \mathrm{weakly\;\; in\;} L^2(L_\xi;\R^{2 \times 2}),\\
&\curdn \tilde{B}_n \wk{} \alpha \quad \mathrm{weakly\;\; in\;}  L^2(L_\xi;\R^{2}),
\end{align}
for some $B \in L^2(L_\xi;\R^{2 \times 2})$ and $\alpha \in L^2(L_\xi;\R^{2})$.

Furthermore, using the quadratic bounds on $\tilde{k}$ and \eqref{eq:k_b_relation}, we deduce that
\begin{align}
	&\tilde{k}_n \wk{} \tilde{k} \quad \mathrm{weakly\;\; in\;}  L^2(L_\xi;\R^{2}), \label{eq:wk_k_tilde}\\
	&\e_n \dgradn \tilde{k}_n \wk{} B \quad  \mathrm{weakly\;\; in\;} L^2(L_\xi;\R^{2 \times 2}), \label{eq:wk_dgrad}
\end{align}
for some $\tilde{k} \in  L^2(L_\xi;\R^{2})$. 

In order to further characterize $B$, we can analyze component-wise for $i = 1,2$ using \eqref{eq:wk_dgrad}. To be precise, for $\phi \in C^\infty_c(L_\xi)$, it follows that
\begin{align*}
    &\int_{L_\xi} B^{i1} \phi\; dx= \limn\int_{L_\xi} \e_n \partial_1 \tilde{k}_n^{i} \phi \;dx = - \limn\e_n \int_{L_\xi} \tilde{k}^i_n \partial_1\phi\;dx = 0, \\
    &\int_{L_\xi} B^{i2} \phi\; dx= \limn\int_{L_\xi} \partial_2 \tilde{k}_n^{i} \phi \;dx = - \limn \int_{L_\xi} \tilde{k}^i_n \partial_2\phi\;dx = \int_{L_\xi} \partial_2 \tilde{k}^i \phi\;dx,
\end{align*}
where we have applied \eqref{eq:wk_k_tilde} after integrating by parts.

Thus, we conclude that 
$$B = \begin{bmatrix}
	0 & \partial_2 \tilde{k}_1 \\
	0 & \partial_2 \tilde{k}_2
\end{bmatrix}.$$ 

Now we can get information on $\curdn \tilde{B}_n$ by integrating by parts. We will do it component-wise for $i = 1,2$. For any $\phi \in C^\infty(L_\xi)$ not necessarily compactly supported, we have
\begin{align*}
	&\int_{L_\xi}{\alpha^i \phi\; dx} = \limn \int_{L_\xi}{[\curdn \tilde{B}_n]^i\;\phi\; dx},\\
	 & = \limn - \int_{L_\xi}{\left[\tilde{B}^{i2}_n\partial_1\phi-\frac{1}{\e_n}\tilde{B}^{i1}_n\partial_2 \phi\right]\; dx} + \int_{\partial L_\xi} \phi [\tilde{B}^{i2}_n\nu_1 -\frac{1}{\e_n}\tilde{B}^{i1}_n\nu_2 ]\; d\mathcal{H}^1,
\end{align*}
where $\mathcal{H}^1$ is the one dimensional Hausdorff (surface) measure in $\R^2$.

Using the tangential relations in \eqref{eq:B_tan} and the weak convergence of $\tilde{B}_n$, we simplify
\begin{multline}\label{eqn:smooth_int_parts_alpha}
	\int_{L_\xi}{\alpha^i \phi\; dx} + \int_{L_\xi}{\partial_2 \tilde{k}^i \partial_1 \phi}\;dx \\= \limn \left[\int_{L_\xi}{\frac{1}{\e_n}\tilde{B}_{i1}^n\partial_2 \phi\; dx} + \int_{-\frac{\xi}{2}}^{\frac{\xi}{2}} \phi(1,x_2) \tilde{B}^{i2}_n(1,x_2) dx_2\right]. 
\end{multline}
Taking $\phi \equiv 1$ leads to the relation
\begin{equation}\label{eqn:alpha_avg}
	\int_{L_\xi}{\alpha^i \; dx} = \limn  \int_{-\frac{\xi}{2}}^{\frac{\xi}{2}} \tilde{B}^{i2}_n(1,x_2) dx_2 .
\end{equation}
Allowing $\phi\in C^\infty((-\xi,1))$, which means that $\partial_2 \phi = 0$, leads to the equation
\begin{align}\label{eqn:partial20}
	\int_{L_\xi}{\alpha^i \phi\; dx} + \int_{L_\xi}{\partial_2 \tilde{k}^i \partial_1 \phi\;dx}= \limn  \phi(1)\int_{-\frac{\xi}{2}}^{\frac{\xi}{2}}\tilde{B}^{i2}_n(1,x_2) dx_2 .
\end{align}
Define $[\tilde{k}^i]: (-\xi,1) \to \R$ as $$[\tilde{k}^i](x_1) := \int_{-\frac{\xi}{2}}^{\frac{\xi}{2}}{\partial_2\tilde{k}^i\;dx_2}.$$
This is well-defined as an $L^2$ function since $\partial_2 \tilde{k}^i \in L^2(L_\xi)$. 
Using \eqref{eqn:alpha_avg} and the fact that $\phi$ only depends on $x_1$, we can simplify the relation \eqref{eqn:partial20} further as
\begin{align}\label{eqn:partial21}
	\int_{L_\xi}{\alpha^i \phi\; dx} + \int_{-\xi}^1{[\tilde{k}^i]\partial_1 \phi\;dx_1}=  \phi(1)\int_{L_\xi}{\alpha^i \; dx} .
\end{align}
In particular, taking $\phi \in C_c^\infty((-\xi,1))$, we deduce that 
$[\tilde{k}^i] \in W^{1,2}((-\xi,1))$ with
$$\frac{d}{dx_1}[\tilde{k}^i](x_1) = \int_{-\frac{\xi}{2}}^{\frac{\xi}{2}} \alpha^i \;dx_2.$$
Now, for generic $\phi\in C^\infty((-\xi,1))$, we can integrate by parts in \eqref{eqn:partial21} to get
\begin{align}
	&[\tilde{k}^i](1)\phi(1)-[\tilde{k}^i](-\xi)\phi(-\xi) = \phi(1)([\tilde{k}^i](1)-[\tilde{k}^i](-\xi)),\\
	&\text{and so  }\qquad [\tilde{k}^i](-\xi)(\phi(1)-\phi(-\xi)) = 0.
\end{align}
Since the equation has to hold for every such $\phi$, we have that $[\tilde{k}^i](-\xi) = 0.$
This gives us a complete characterization of the vertical jump as
\begin{equation}
	\label{eqn:k_jump}
	[\tilde{k}^i](s) = \int_{-\xi}^{s}\int_{-\frac{\xi}{2}}^{\frac{\xi}{2}} \alpha^i \;dx.
\end{equation}

By our convergences, we also have that \eqref{eq:curl_contraint_cov} passes to the limit. To be precise,
\begin{equation}\label{eq:curl_jump}
	2 = \left|\int_{L_\xi}{\curd \tilde{B}\;dx} \right| \to \left|\int_{L_\xi}{\alpha \;dx} \right| = \left|[\tilde{k}](1)\right|  .
\end{equation}
\end{proof}

\subsection{Coincidence of traces}
We show that the compatibility relation from Theorem \ref{thm:compactness} relates the jump of $k$ directly to the limit defect field $\alpha.$
\begin{proposition}\label{prop:traceCoincidence}
Supposing that $k,\tilde k$ arise as in Theorem \ref{thm:compactness}, then the compatibility relation
$$\jump{k} = [\tilde{k}]  \quad \text{ for } \mathcal{H}^1\text{-almost every }x\in (-\xi,1)\times \{0\}$$
is satisfied, where $\jump{k}$ denotes the $BV$ jump, oriented as the trace from $\{x_2>0\}$ minus the trace from $\{x_2<0\}$.
\end{proposition}

\begin{proof}

We begin by noting that the function $k$ is such that $$Dk = \nabla k \mathcal{L}^2 + \jump{k}\otimes e_2 d\mathcal{H}^{1}\llcorner_{(-\xi,1)\times\{0\}}.$$
From this it follows that for $i=1,2$ and $\phi \in C_c^\infty (\Omega)$, we have
\begin{equation}\nonumber
\int_{\Omega} \partial_2 \phi k^i\, dx = - \left(\int_\Omega \phi \partial_2 k^i \, dx + \int_{(-\xi,1)\times\{0\}} \phi\jump{k^i} \, \mathcal{H}^{1} \right).
\end{equation}
Consequently, to prove the theorem, it suffices to show that for all $\phi \in C_c^\infty (\Omega \setminus \{(-\xi,0)\})$ it holds that
\begin{equation}\label{eqn:derivativeClaim}
\int_{\Omega} \partial_2 \phi k^i\, dx = - \left(\int_\Omega \phi \partial_2 k^i \, dx + \int_{(-\xi,1)\times\{0\}} \phi [\tilde{k}]^i \, \mathcal{H}^{1} \right).
\end{equation}

Let $\{k_n\}$ be the sequence from Theorem \ref{thm:compactness}. From here, we drop the superscript $i$ but continue to operate component-wise. Directly by the uniform $L^2(\Omega;\R^2)$ bound on $k_n$ and strong convergence away from $(-\xi,1)\times\{0\},$ we have that
\begin{equation}\label{eqn:limitBeforeIBP}
\int_{\Omega} \partial_2 \phi k_n\, dx \to \int_{\Omega} \partial_2 \phi k\, dx \quad\text{ as }n \to \infty.
\end{equation}
Performing an integration by parts, we also find
\begin{equation}\label{eqn:integralSplit}
\int_{\Omega} \partial_2 \phi k_n\, dx =  - \int_{\Omega} \phi  \partial_2 k_n\, dx = - \left(\int_{\Omega \setminus L_{\e_n,\xi}} \phi  \partial_2 k_n\, dx + \int_{L_{\e_n,\xi}}\phi  \partial_2 k_n\, dx\right).
\end{equation}
By (\ref{eq:Uniform_bds_Bulk}), $\nabla k_n\chi_{\Omega\setminus L_{\e_n,\xi}} \wkto \nabla k$ in $L^2(\Omega).$ The first term on right-hand side of (\ref{eqn:integralSplit}) converges with 
\begin{equation}\label{eqn:bulkLimit}
\int_{\Omega \setminus L_{\e_n,\xi}} \phi  \partial_2 k_n\, dx \to \int_{\Omega } \phi  \partial_2 k\, dx \quad \text{ as }n \to \infty.
\end{equation}

For the second term, we recall that $\tilde{k}_n (y_1,y_2) = k_n(y_1,\e y_2),$ and perform a change of variables (recall $L_{\xi} : = L_{1,\xi}$)
\begin{equation}\label{eqn:layerBlowup}
\begin{aligned}
\int_{L_{\e_n,\xi}}\phi  \partial_2 k_n\, dx =& \int_{L_{\xi}}\phi(y_1,\e y_2) \partial_2 k_n(y_1,\e y_2) \e \, d(y_1,y_2)\\
=&  \int_{L_{\xi}}\phi(y_1,\e y_2) \partial_2 \tilde{k}_n(y)  \, d(y_1,y_2).\\
\end{aligned}
\end{equation}
Given the regularity of $\phi, $ the map $\phi_\e(y) : = \phi(y_1,\e y_2)$ converges strongly in $L^2(L_{\xi})$ to $\phi(y_1,0)$. By (\ref{eq:k_b_relation}) and (\ref{eq:B_bounds}), it follows that $\partial_2 \tilde{k}_n \wkto \partial_2 \tilde{k}$ in $L^2(L_{\xi})$. Passing to the limit in (\ref{eqn:layerBlowup}) as $n \to \infty$, applying Fubini's theorem, and recalling the definition (\ref{def:tildekjump}) for $[\tilde k]$, we find 
\begin{equation}\label{eqn:layerLimit}
\int_{L_{\e_n,\xi}}\phi  \partial_2 k_n\, dx \to \int_{L_{\xi}}\phi(y_1,0)  \partial_2 \tilde{k}\, dy = \int_{(-\xi,1)\times\{0\}} \phi [\tilde{k}] \, d\mathcal{H}^{1}.
\end{equation}

Finally, to obtain (\ref{eqn:derivativeClaim}), we pass to the limit in (\ref{eqn:integralSplit}) using (\ref{eqn:limitBeforeIBP}), (\ref{eqn:bulkLimit}), and (\ref{eqn:layerLimit}), thereby concluding the proposition.
\end{proof}
\begin{corollary}\label{Cor: CurlNablak}
    Letting $k$ be as in Theorem \ref{thm:compactness}, $\cur {\nabla k}$ is a finite Radon measure.
\end{corollary}
\begin{proof}
We treat $\cur \nabla k$ row-wise.

We start from the relation $\cur(Dk^i)=0$ in the sense of distributions. Using the decomposition of derivatives for SBV functions, we can write for any $\phi \in C^\infty_c(\o)$:
\begin{equation}
    \langle{\cur(\nabla k^i \mathcal{L}^2)},\phi\rangle = - \langle{\cur(\jump{k^i} e_2 \mathcal{H}^{1}\llcorner_{(-\xi,1)\times\{0\}})},\phi\rangle.
\end{equation}
where the brackets denote the duality pairing of distributions.
On the left hand side, by definition, 
$$\langle{\cur(\nabla k^i \mathcal{L}^2)},\phi\rangle = \langle\cur\nabla k^i,\phi\rangle.$$
On the right hand side, we also unwrap the duality and use the area formula to find
\begin{equation}\label{eqn:jump_curl}
    - \langle{\cur(\jump{k^i} e_2 \mathcal{H}^{1}\llcorner_{(-\xi,1)\times\{0\}})},\phi\rangle = \int_{-\xi}^1\jump{k^i}(x_1)\partial_1\phi(x_1,0)\;dx_1.
\end{equation}
But by Theorem \ref{thm:compactness} and Proposition \ref{prop:traceCoincidence}, we know that $\jump{k^i} \in W^{1,2}((-\xi,1)$ so we can integrate by parts and achieve
\begin{align}
    \left|\langle{\cur(\jump{k^i} e_2 \mathcal{H}^{1}\llcorner_{(-\xi,1)\times\{0\}})},\phi\rangle\right|&= \left|\int_{-\xi}^1\frac{d}{dx_1}\jump{k^i}(x_1)\phi(x_1,0)\;dx_1 \right|\nonumber\\
    &\leq \left\|\frac{d}{dx_1}\jump{k^i}\right\|_{L^2(-\xi,1)}\|\phi\|_{\infty}\leq C\|\phi\|_{\infty}.
\end{align}
Thus combining the previous equations, we conclude that
$$\left |\langle\cur\nabla k^i,\phi\rangle\right |\leq C\|\phi\|_{\infty}.$$
Therefore, $\cur\nabla k^i$ is a finite Radon measure, and furthermore the integration by parts argument shows that $\cur\nabla k^i = - \frac{d}{dx_1}\jump{k^i}e_2\mathcal{H}^{1}\llcorner_{(-\xi,1)\times\{0\}}.$
\end{proof}

\begin{proof}[Proof of Theorem \ref{thm:compactCorollary}]
This follows directly from Theorem \ref{thm:compactness}, Proposition \ref{prop:traceCoincidence}, and Corollary \ref{Cor: CurlNablak}.
\end{proof}

\section{Conjectures about the Limiting Energy}\label{sec:conclusion}
Though Theorems \ref{thm:relax} and \ref{thm:compactCorollary} address the behavior of the energy $E_{\e,\xi}$, they leave the $\Gamma$-limits of in the singularly perturbed regimes ($\e,\xi \to 0$) unresolved. In the constrained setting, the principal challenge to characterize the limiting energy is to understand the coupled term $$\int_{L_\xi} \e_n\left|\nabla_{\e_n} \tilde{k}_n - \frac{\tilde{B}_n}{\e_n}\right|^2 \;dx.$$

If we rewrite this term in a rescaled Helmholtz decomposition \\ $\tilde{B}_n = \dgradn z_n + p_n$, where $p_n$ will be rescaled divergence free, using orthogonality of $p_n$ with respect to $\nabla_{\e_n}z_n$, this coupled term will become
$$\int_{L_\xi} \left[\e_n\left|\nabla_{\e_n} \tilde{k}_n - \frac{\dgradn z_n}{\e_n}\right|^2 + \frac{|p_n|^2}{\e_n} \right] dx.$$
So we see that if this has bounded energy, then we will have $p_n \to 0$ strongly and the limit of $\tilde{B}_n$ comes purely from the rescaled gradient term.

Using this, one can heuristically argue for the structure of the limiting energy as follows. Considering the terms which depend on only $\tilde{B}_n$, we can view the relaxation in the double well-function as similar to the dimension reduction where we fix the so called bending vector $\frac{1}{\e}\partial_2 z_n \wk{} \partial_2 \tilde{k}$. Such a relaxation has been considered in the 3D-2D case in \cite{BouFonMasc}, and gives a cross convex-quasiconvex envelope (see also \cite{FonKindPed}). Using the convergences given in Theorem \ref{thm:compactness}, one could imagine leveraging the limiting structure of $B$ and weak lower semicontinuity of cross convex-quasiconvex envelopes (denoted here by $Q^*$) with respect to rescaled gradients
\begin{align*}
    &\int_{L_\xi}\left[\xi^2|\cur_{\e_n}\tilde{B}_n|^2 + \frac{1}{\xi^2}W(\xi |\tilde{B}_n|)\right]\; dx\\
    &\to \int_{L_\xi}\left[\xi^2|\alpha|^2 + \frac{1}{\xi^2}Q^*(W(|\cdot|))(\xi \partial_2\tilde{k}\otimes e_2)\right]\; dx.
\end{align*}
Since the envelope \eqref{eq:limit} we are considering should not depend on $\alpha,\tilde{k}$, we may optimize with respect to a lower bound achieved. In particular recalling definition \ref{def:tildekjump} and Proposition \ref{prop:traceCoincidence}, we see through Jensen's inequality and the definition of cross quasiconvexity a possible lower bound is
\begin{multline}\label{conj:jump}
    \int_{L_\xi}\left[\xi^2|\alpha|^2 + \frac{1}{\xi^2}Q(W(|\cdot|))(\xi \partial_2\tilde{k}\otimes e_2)\right]\; dx \\ 
    \geq \int_{-\xi}^1{\left[\xi\left|\frac{d}{dx_1}\jump{k}(x_1)\right|^2 + \frac{1}{\xi}Q(W(|\cdot|))(\jump{k})\right]}\; dx.
\end{multline}

In this setting, this equation is quite similar to a Modica-Mortola functional for the vertical jump of the director with a transition layer on the order of the core length $\xi$. This is the type of picture predicted by the numerical experiments in \cite{ZANV}. 

However, if the correct energy is as in (\ref{conj:jump}), recalling Remark \ref{rmk:quasiconvex}, we see that the disclination layer $(-\xi,1)\times\{0\}$ allows for any jump of the unit-length director field. To counter this, most likely the energy should be modified so to obtain strong convergence (in at least $L^1$) of $\tilde{B}_n$, so that the quasiconvexification of $W(|\cdot|)$ doesn't occur. In principle, this could be done via the inclusion of a smaller order term $\int_{\Omega} \e^3|{\rm div} B|^2 \, dx$ within the energy. Given the smaller order of $\e$, it would most likely be negligible in the limit $\e \to 0$. However, at the $\e>0$ fixed level, it will require that $B=0$ in the sense of traces on the boundary of the layer. This has the effect that $\nabla k$ will still see some small energy contribution within the layer. 

\section*{Acknowledgment}
This work was supported by the grant NSF OIA-DMR \#2021019.
L.G was partially supported by the Deutsche Forschungsgemeinschaft 320021702 GRK2326 Energy, Entropy, and Dissipative Dynamics (EDDy). K.S. also acknowledges funding by the Deutsche Forschungsgemeinschaft (DFG, German Research Foundation) under Germany's Excellence Strategy - GZ 2047/1, Projekt-ID 390685813.
\bibliographystyle{alphaurl}\bibliography{layer_math.bib}
\end{document}